\documentclass{article}
\usepackage{amsmath,amsthm} 

\def\today{June 26, 2001}

\advance \topmargin by -\headheight
\advance \topmargin by -\headsep
     
\evensidemargin \oddsidemargin
\theoremstyle{plain}  

\newtheorem{theorem}{Theorem}
\newtheorem{lemma}[theorem]{Lemma}
\newtheorem{e-proposition}[theorem]{Proposition}

\newtheorem{theoreme}{Th\'eor\`eme}

\theoremstyle{definition}

\newtheorem{e-definition}[theorem]{Definition}

\theoremstyle{remark}

\newtheorem{remark}[theorem]{Remark}

\newcommand{\PU}{\mathrm{PU}}
\newcommand{\U}{\mathrm{U}}
\newcommand{\SU}{\mathrm{SU}}
\DeclareMathOperator{\ad}{ad}
\DeclareMathOperator{\rk}{rk}
\DeclareMathOperator{\Hom}{Hom}
\DeclareMathOperator{\End}{End}

\begin{document}

\vspace*{1cm}
\noindent
{\LARGE \textbf{Representations of the fundamental group of a surface
    in $\PU(p,q)$ and holomorphic triples} \par}

\vspace{5mm}
\noindent
{\bfseries
Steven B. Bradlow~$^{\text{a}}$,\ \
Oscar Garc{\'\i}a-Prada~$^{\text{b}}$,\ \
Peter B. Gothen~$^{\text{c}}$
}

\vspace{5mm}
\begin{itemize}\labelsep=2mm\leftskip=-5mm
{\bfseries
\item[$^{\text{a}}$]
Department of Mathematics,
University of Illinois,
Urbana,
IL 61801,
USA \\
E-mail: \texttt{bradlow@math.uiuc.edu}
\item[$^{\text{b}}$]
Departamento de Matem{\'a}ticas, 
Universidad Aut{\'o}noma de Madrid,
28049 Madrid,
Spain\\
E-mail: \texttt{oscar.garcia-prada@uam.es}
\item[$^{\text{c}}$]
Departamento de Matem{\'a}tica Pura, 
Faculdade de Ci{\^e}ncias,
4099-002 Porto, 
Portugal\\
E-mail: \texttt{pbgothen@fc.up.pt}
}
\end{itemize}

\noindent
{\bfseries(\today)}

\noindent
\hrulefill


\noindent
\textbf{Abstract.}  We count the connected components in the moduli
space of $\PU(p,q)$-representations of the fundamental group for a
closed oriented surface. The components are labelled by pairs of
integers which arise as topological invariants of the flat bundles
associated to the representations. Our results show that for each
allowed value of these invariants, which are bounded by a Milnor--Wood
type inequality, there is a unique non-empty connected component.
Interpreting the moduli space of representations as a moduli space of
Higgs bundles, we take a Morse theoretic approach using a certain
smooth proper function on the Higgs moduli space. A key step is the
identification of the function's local minima as moduli spaces of
holomorphic triples. We prove that these moduli spaces of triples are
non-empty and irreducible.

\vspace{5mm}

{\itshape {\bfseries\large \noindent Repr\'esentations du groupe
    fondamental d'une surface dans $\PU(p,q)$ et tripl\'es holomorphes
    \par }

\vspace{2mm}


\noindent
{\bfseries{}R{\'e}sum{\'e.}}} Nous comptons le nombre de composantes
connexes de l'espace des modules de repr\'esentations du groupe
fondamental d'une surface compacte orient\'ee dans $\PU(p,q)$.  Les
invariants topologiques du fibr\'e plat correspondant permettent
d'associer \`a une telle repr\'esentation deux entiers, dont la valeur
est born\'ee par une in\'egalit\'e du type Milnor--Wood.  Nous
montrons que pour n'importe quelle valeur autoris\'ee de ces
invariants, il existe une unique composante connexe non vide dans
l'espace des modules. Notre m\'ethode utilise l'interpr\'etation de
l'espace des modules de repr\'esentations comme un espace des modules
de fibr\'es de Higgs. L'analyse des minimum locaux d'une certaine
fonction propre d\'efinie sur l'espace des modules est un \'el\'ement
cl\'e de la d\'emonstration~: nous identifions ces ensembles avec
l'espace des modules de tripl\'es holomorphes et d\'emontrons que ces
espaces de modules sont non vides et irr\'eductibles.

\noindent
\hrulefill


\vspace{5mm}

\noindent
{\bfseries\itshape\large{}Version fran\c{c}aise abr\'eg\'ee}
\vspace{2mm}

Soit $X$ une surface compacte de genre $g \geq 2$ et soit $G$ un
groupe de Lie connexe.  L'ensemble des repr\'esentations r\'eductives
$\Hom^{+}(\pi_1 X, G)$ du groupe fondamental de $X$ dans $G$ est une
vari\'et\'e analytique r\'eelle et par cons\'equent un espace
topologique avec une topologie analytique. On a une action par
conjugaison de $G$ sur $\Hom^{+}(\pi_1 X, G)$ et on munit l'espace des
modules
\begin{displaymath}
  \mathcal{M}_{G} = \Hom^{+}(\pi_1 X, G) / G
\end{displaymath}
de la topologie quotient; c'est un espace de Hausdorff parce que nous
ne consid\'erons que des repr\'esentations r\'eductives.  Il est bien
connu que $\mathcal{M}_{G}$ peut \'egalement \^etre interpr\'et\'e
comme l'espace des modules de $G$-fibr\'es plats sur $X$.  Cet espace
qui est une source d'int\'er\^et tant en math\'ematiques qu'en
physique a \'et\'e \'etudi\'e largement du point de vue de la
topologie, des th\'eories de jauge et de la g\'eom\'etrie
alg\'ebrique.  L'objet de cet article est de d\'eterminer le nombre de
composantes connexes de $\mathcal{M}_{G}$ dans le cas o\`u $G =
\PU(p,q)$.  Des r\'esultats ant\'erieurs concernant ce probl\`eme pour
des groupes r\'eels non compacts $G$ sont dus \`a Goldman
\cite{goldman:1988}, Gothen \cite{gothen:2001}, Hitchin
\cite{hitchin:1992}, Markman et Xia \cite{markman-xia:2000}, et Xia
\cite{xia:1997,xia:1999,xia:2000}.

Afin d'\'enoncer notre r\'esultat, nous devons faire quelques
remarques pr\'eliminaires.  Soit $[\rho] \in \mathcal{M}_{\PU(p,q)}$,
alors le $\PU(p,q)$-fibr\'e plat correspondant peut \^etre relev\'e en
un $\U(p,q)$-fibr\'e qui admet de plus une r\'eduction \`a son sous
groupe compact maximal $\U(p) \times \U(q)$.  Soit $V \oplus W$ le
fibr\'e vectoriel complexe de rang $p+q$ correspondant, alors sa
classe d'isomorphisme est ind\'ependante des choix que nous avons
faits et le couple d'entier $(d_V,d_W) = \bigl(\deg(V),\deg(W)\bigr)$
est un invariant $[\rho]$ bien d\'efini.  D'apr\`es les r\'esultats de
Domic et Toledo \cite{domic-toledo:1987} (am\'eliorant une borne
obtenue par Dupont \cite{dupont:1978} dans le cas o\`u $G =
\SU(p,q)$), on a l'in\'egalit\'e du type Milnor-Wood
\begin{equation}
\label{eq:f-milnor-wood}
  \frac{\lvert q d_V - p d_W \rvert}{p+q}
  \leq \min\{p,q\}(g-1), 
\end{equation}
donnant des bornes sur les valeurs 
que l'invariant topologique 
$(d_V,d_W)$ peut prendre. D\'efinissons
\begin{displaymath}
  \mathcal{M}(d_V,d_W) =
  \{[\rho] \; : \; \bigl(\deg(V),\deg(W)\bigr) = (d_V,d_W) \},
\end{displaymath}
o\`u $V \oplus W$ est le fibr\'e vectoriel de rang  $p+q$  associ\'e \`a
$[\rho]$ comme ci-dessus.  Nous pouvons maintenant \'enoncer notre
r\'esultat principal.
\begin{theoreme}
  Pour chaque  $(d_V,d_W) \in \mathbf{Z}^{2}$ satisfaisant 
  \eqref{eq:f-milnor-wood} 
 l'espace $\mathcal{M}(d_V,d_W)$ est non vide
  et c'est une composante connexe de l'espace des modules
  $\mathcal{M}_{\PU(p,q)}$.
\end{theoreme}

Il est plus naturel si on adopte le point de vue des repr\'esentations
du groupe fondamental et des fibr\'e plats de travailler avec la forme
associ\'ee $\PU(p,q)$ de $\U(p,q)$.  Pourtant on peut tout de m\^eme
consid\'erer l'espace des modules $\mathcal{M}_{\U(p,q)}$ des fibr\'es
munis d'une connexion r\'eductive avec courbure centrale constante; en
terme de repr\'esentations du groupe fondamental ceci correspond aux
repr\'esentations d'une extension centrale universelle de $\pi_1 X$
dans $\U(p,q)$. On a une application naturelle $\mathcal{M}_{\U(p,q)}
\to \mathcal{M}_{\PU(p,q)}$ dont les fibres s'identifient \`a
$\Hom(\pi_1 X, \U(1))$.  Puisque le centre $\U(1)$ de $\U(p,q)$ est
connexe, cette application induit un isomorphisme entre les ensembles
de composantes connexes des espaces de modules.  Par cons\'equent,
nous pourrons \'egalement consid\'erer l'espace
$\mathcal{M}_{\U(p,q)}$.  Ceci est finalement plus naturel du point de
vue des techniques sur les fibr\'es holomorphes que nous emploierons.
La d\'emostration repose sur un approche holomorphe et n\'ecessite de
donner \`a $X$ une structure de surface de Riemann. La th\'eorie de
fibr\'es de Higgs et leurs liens avec les fibr\'es plats (Corlette
\cite{corlette:1988}, Donaldson \cite{donaldson:1987}, Hitchin
\cite{hitchin:1987, hitchin:1992}, Simpson
\cite{simpson:1988,simpson:1992,simpson:1994a,simpson:1994b}) ansi que
les tripl\'es holomorphes (\cite{bradlow-garcia-prada:1996},
\cite{garcia-prada:1994}) interviennent de fa\c{c}on essentielle.

\par\medskip\centerline{\rule{2cm}{0.2mm}}\medskip
\setcounter{section}{0}
\setcounter{equation}{0}


\subsubsection*{1. Introduction}

Let $X$ be a closed oriented surface of genus $g \geq 2$ and let $G$
be a connected Lie group.  The set of reductive representations
$\Hom^{+}(\pi_1 X, G)$ of the fundamental group of $X$ in $G$ is a
real analytic variety, and hence it is a topological space with the
analytic topology.  There is an action of $G$ on $\Hom^{+}(\pi_1 X,
G)$ by conjugation and we give the moduli space
\begin{displaymath}
  \mathcal{M}_{G} = \Hom^{+}(\pi_1 X, G) / G
\end{displaymath}
the quotient topology; this is Hausdorff because we have restricted
attention to reductive representations.  It is well known that
$\mathcal{M}_{G}$ can also be viewed as the moduli space of flat
$G$-bundles on $X$.  This space has been the object of much interest
both in mathematics and physics and it has been extensively studied
from the points of view of topology, gauge theory and algebraic
geometry.
The purpose of the present note is to determine the number of
connected components of $\mathcal{M}_{G}$ in the case of $G =
\PU(p,q)$.  Previous results on this problem for non-compact real
groups $G$ include those of Goldman \cite{goldman:1988},
Gothen \cite{gothen:2001}, Hitchin \cite{hitchin:1992},
Markman and Xia \cite{markman-xia:2000}, and Xia
\cite{xia:1997,xia:1999,xia:2000}.

In order to state our result we need some preliminary observations.
Let $[\rho] \in \mathcal{M}_{\PU(p,q)}$, then the corresponding flat
$\PU(p,q)$-bundle can be lifted to a $\U(p,q)$-bundle which, in turn,
has a reduction to its maximal compact subgroup $\U(p) \times \U(q)$.
Let $V \oplus W$ be the corresponding complex rank $p+q$ vector
bundle, then its isomorphism class is independent of the choices made
and hence the pair of integers $(d_V,d_W) =
\bigl(\deg(V),\deg(W)\bigr)$ is a well defined invariant of $[\rho]$.
It follows from results of Domic and Toledo \cite{domic-toledo:1987}
(improving on a bound obtained by Dupont \cite{dupont:1978} in the
case $G = \SU(p,q)$), that there is the Milnor-Wood type inequality
\begin{equation}
\label{eq:milnor-wood}
  \frac{\lvert q d_V - p d_W \rvert}{p+q}
  \leq \min\{p,q\}(g-1),
\end{equation}
giving bounds on the possible values of the topological invariants
$(d_V,d_W)$.  Define
\begin{displaymath}
  \mathcal{M}(d_V,d_W) =
  \{[\rho] \; : \; \bigl(\deg(V),\deg(W)\bigr) = (d_V,d_W) \},
\end{displaymath}
where $V \oplus W$ is the rank $p+q$ vector bundle associated to
$[\rho]$ as above.  We can now state our main result.
\begin{theorem}
  \label{thm:main}
  For each $(d_V,d_W) \in \mathbf{Z}^{2}$ satisfying
  \eqref{eq:milnor-wood} the space $\mathcal{M}(d_V,d_W)$ is non-empty
  and it is a connected component of the moduli space
  $\mathcal{M}_{\PU(p,q)}$.
\end{theorem}

From the points of view of representations of the fundamental group
and flat bundles it is most natural to work with the adjoint form
$\PU(p,q)$ of $\U(p,q)$.  However one can also consider the moduli
space $\mathcal{M}_{\U(p,q)}$ of bundles with a reductive connection
with constant central curvature; in terms of representations of the
fundamental group this corresponds to representations of a universal
central extension of $\pi_1 X$ in $\U(p,q)$.  There is a natural map
$\mathcal{M}_{\U(p,q)} \to \mathcal{M}_{\PU(p,q)}$ whose fibres can be
identified with $\Hom(\pi_1 X, \U(1))$.  Since the centre $\U(1)$ of
$\U(p,q)$ is connected, this map induces an isomorphism between the
sets of connected components of the moduli spaces.  Thus we could
equally well consider the space $\mathcal{M}_{\U(p,q)}$.  This is
actually more natural from the point of view of the holomorphic vector
bundle methods which we shall employ, and we shall do this from now
on.

\subsubsection*{2. Outline of proof of Theorem \ref{thm:main}}

The proof rests on a holomorphic approach.  This requires giving $X$
the structure of a Riemann surface.  The two basic ingredients are
Higgs bundles and holomorphic triples.

A \emph{Higgs bundle} is a pair $(E,\Phi)$, where $E$ is a holomorphic
bundle on $X$ and $\Phi \in H^0(X; \End(E) \otimes K)$ is a
holomorphic endomorphism of $E$ twisted by the canonical bundle $K$ of
$X$.  The slope of a Higgs bundle is by definition the slope of the
underlying vector bundle, $\mu(E) = \deg(E)/\rk(E)$.  There is a
stability condition, generalizing the usual slope stability for
holomorphic vector bundles: a Higgs bundle is called stable if $\mu(F)
< \mu(E)$ for any proper subbundle $F$ of $E$ preserved by $\Phi$, and
one has the corresponding notions of semi-stability and
poly-stability.

By results of Corlette \cite{corlette:1988}, Donaldson
\cite{donaldson:1987}, Hitchin \cite{hitchin:1987, hitchin:1992}, and
Simpson \cite{simpson:1988,simpson:1992,simpson:1994a,simpson:1994b}
there is a homeomorphism between the moduli space
$\mathcal{M}_{\U(p,q)}$ and the moduli space of poly-stable Higgs
bundles $(E,\Phi)$ of the form
\begin{displaymath}
  \label{eq:upq-higgs-bundle}
  E = V \oplus W, \qquad
  \Phi =
  \left(
  \begin{smallmatrix}
    0 & b \\
    c & 0
  \end{smallmatrix}
  \right)
\end{displaymath}
where $V$ and $W$ are holomorphic vector bundles on $X$ of rank
$p$ and $q$ respectively, $b \in
H^0(\Hom(W,V) \otimes K)$, and $c \in H^0(\Hom(V,W) \otimes
K)$.  We shall call such a Higgs bundle a \emph{$\U(p,q)$-Higgs
  bundle}.

The homeomorphism above comes via an interpretation of
$\mathcal{M}_{\U(p,q)}$ as a gauge theory moduli space: it is possible
to find a preferred (harmonic) metric in a reductive flat bundle and
also to find a preferred metric in a poly-stable Higgs bundle.  Both
of these constructions lead to a solution to a set of gauge theoretic
equations known as Hitchin's equations \cite{hitchin:1987}.  Thus
$\mathcal{M}_{\U(p,q)}$ can also be seen as the gauge theory moduli
space of solutions to Hitchin's equations modulo gauge equivalence.
This point of view allows one to define the function
\begin{displaymath}
  f(E,\Phi) = \int_{X} \lvert \Phi \rvert^{2}
\end{displaymath}
on $\mathcal{M}_{\U(p,q)}$, which from Uhlenbeck's weak compactness
theorem is known to be proper \cite{hitchin:1987}.

Clearly each $\mathcal{M}(d_V,d_W)$ is a union of connected
components.  Hence all we need to prove is that each of these spaces
is non-empty and connected.  Since $f$ is proper it is sufficient to
show connectedness of the subspace $\mathcal{N}(d_V,d_W)$ of local
minima of $f$ restricted to $\mathcal{M}(d_V,d_W)$.  Thus the first
main step in the proof of Theorem \ref{thm:main} is the following
characterization of the local minima of $f$.

\begin{lemma}
    \label{lemma:minima}
  A poly-stable $\U(p,q)$-Higgs bundle is a local minimum of
  $f$ if and only if at least one of the sections $b$ and $c$
  vanishes.
\end{lemma}

\begin{proof}[Proof (sketch)]
  Assume that $(E,\Phi)$ is a stable $\U(p,q)$-Higgs bundle.  This
  means that it represents a smooth point of the moduli space (an
  extra argument, similar to the one given by Hitchin in
  \cite{hitchin:1992}, is required to deal with the case of non-smooth
  points of the moduli space).  It is known from \cite{hitchin:1992}
  that $(E,\Phi)$ represents a critical point of $f$ on
  $\mathcal{M}_{\U(p,q)}$ if and only if it is a variation of Hodge
  structure.  In our case this means that it is a $\U(p,q)$-Higgs
  bundle of the form
  \begin{displaymath}
    E = F_1 \oplus \cdots \oplus F_m,
  \end{displaymath}
  where $\Phi$ maps $F_k$ to $F_{k+1} \otimes K$ and each $F_{k}$ is
  contained in either $V$ or $W$.
  Define $$U_k =
  \bigoplus_{i-j=k} \Hom(F_j,F_i). $$  It follows from the results of
  \cite{hitchin:1992} that the eigenvalues of the Hessian of 
  $f$ at $(E,\Phi)$ are all even integers and that its $-2k$-eigenspace
  is isomorphic to the first hypercohomology of
  the complex
  \begin{displaymath}
    C^{\bullet}_{2k} : U_{2k}
    \xrightarrow{\ad(\Phi)} U_{2k+1} \otimes K.
  \end{displaymath}
  The key ingredient in the proof is then the following vanishing
  criterion: the first hypercohomology
  $\mathbf{H}^{1}(C^{\bullet}_{2k})$ vanishes if and only if $\ad(\Phi)
  \colon U_{2k} \to U_{2k+1} \otimes K$ is an
  isomorphism.  This is proved using the fact that
  $\bigl(\End(E),\ad(\Phi)\bigr)$ is a semi-stable Higgs bundle.
  Using this criterion it is then a matter of elementary linear
  algebra to show that $\mathbf{H}^{1}(C^{\bullet}_{m-1}) \neq 0$.  It
  follows that $(E,\Phi)$ is not a local minimum of $f$ whenever $m
  \geq 3$, thus concluding the proof.
\end{proof}

\begin{remark}
  We can actually be a bit more specific in the statement of this
  Lemma: it is easy to see that a $\U(p,q)$-Higgs bundle with $d_V/p <
  d_W/q$ is a local minimum if and only if $c=0$, that a
  $\U(p,q)$-Higgs bundle with $d_V/p > d_W/q$ is a local minimum if
  and only if $b=0$, and that a $\U(p,q)$-Higgs bundle with $d_V/p =
  d_W/q$ is a local minimum if and only if $b=c=0$.
\end{remark}

\begin{remark}
  The vanishing criterion for hypercohomology described in the above
  proof can also be formulated for groups $G$ other than $\U(p,q)$ and
  it should prove useful for the analysis of the connected components
  of $\mathcal{M}_{G}$ in these cases.
\end{remark}

Lemma \ref{lemma:minima} allows us to reinterpret
$\mathcal{N}(d_V,d_W)$ as a moduli space of so-called holomorphic
triples.  These objects were studied in 
\cite{bradlow-garcia-prada:1996} and \cite{garcia-prada:1994};
we briefly recall the relevant
definitions.
A \emph{holomorphic triple} on $X$,
$T = (E_{1},E_{2},\phi)$ consists of two holomorphic
vector bundles $E_{1}$ and $E_{2}$ on $X$ and a
holomorphic map $\phi \colon E_{2} \to E_{1}$.  
For any $\alpha \in \mathbf{R}$ the 
\emph{$\alpha$-slope} of $T$ is
defined to be
\begin{displaymath}
  \mu_{\alpha}(T)
  = \mu(E_{1} \oplus E_{2}) +
  \alpha\frac{\rk(E_{2})}{\rk(E_{1})+
    \rk(E_{2})}.
\end{displaymath}
There is an obvious notion of subtriple and, using the $\alpha$-slope,
one can define the notions of stability, semi-stability and
poly-stability in the standard way. The existence of moduli spaces of
$\alpha$-stable triples was proved in
\cite{bradlow-garcia-prada:1996}.  Denote by
\begin{displaymath}
  \mathcal{M}_{\alpha}(n_1,n_2,d_1,d_2)
\end{displaymath}
the moduli space of $\alpha$-poly-stable triples $T$ with
$\rk(E_i)=n_i$ and $\deg(E_i) = d_i$ for $i=1,2$.

Let $(E,\Phi)$ be a $\U(p,q)$-Higgs bundle with $c = 0$.  We 
can then define a holomorphic triple 
$(E_{1},E_{2},\phi)$ by setting
\begin{displaymath}
  E_{1} = V\otimes K, \qquad
  E_{2} = W, \qquad
  \phi = b;
\end{displaymath}
and, conversely, given a holomorphic triple we can define an
associated $\U(p,q)$-Higgs bundle with $c = 0$.  Analogously, there is
a bijective correspondence between $\U(p,q)$-Higgs bundles with $b = 0$
and holomorphic triples.  Of course $(n_1,n_2,d_1,d_2)$ can be
expressed in terms of $(p,q,d_V,d_W)$, and vice-versa.

There is a link between the stability conditions for holomorphic
triples and $\U(p,q)$-Higgs bundles: one can show (see 
\cite{gothen:2001}) that a $\U(p,q)$-Higgs bundle $(E,\phi)$ with $b=0$
or $c=0$ is \mbox{(semi-)}stable if and only if the corresponding
holomorphic triple $(E_{1},E_{2},\phi)$ is $\alpha$-(semi-)stable for
$\alpha = 2g-2$.  Lemma \ref{lemma:minima} then implies the following
result.

\begin{lemma}
  \label{lemma:minima-triples}
  The subspace $\mathcal{N}(d_V,d_W)$ of local minima of $f$ on
  $\mathcal{M}(d_V,d_W)$ is isomorphic to the moduli space
  $\mathcal{M}_{\alpha}(n_1,n_2,d_1,d_2)$ of $\alpha$-poly-stable
  triples for $\alpha=2g-2$ (and suitable values of
  $(n_1,n_2,d_1,d_2)$).
\end{lemma}

Thus Theorem \ref{thm:main} follows from our second main result:

\begin{theorem}
  \label{thm:main2}
  The moduli space $\mathcal{M}_{\alpha}(n_1,n_2,d_1,d_2)$ of
  $\alpha$-poly-stable triples is non-empty and irreducible for
  $\alpha \geq 2g-2$ and the values of the topological invariants
  $(n_1,n_2,d_1,d_2)$ allowed by \eqref{eq:milnor-wood}.
\end{theorem}

\begin{proof}[Proof (sketch)]
  The proof of this Theorem applies the strategy used by Thaddeus
  \cite{thaddeus:1994} in his proof of the Verlinde formula.  First
  one obtains a relatively simple description of the moduli space
  $\mathcal{M}_{\alpha}$ for an extreme value of the parameter
  $\alpha$.  Next one studies the variation of the moduli spaces
  $\mathcal{M}_{\alpha}$ as $\alpha$ varies, in order to obtain
  information about $\mathcal{M}_{\alpha}$ for the value of $\alpha$
  in which one is interested.
  
  We note that it is sufficient to consider the case $n_1 \geq n_2$,
  since the case $n_1 \leq n_2$ can be dealt with via duality of
  triples.
  
  Via a careful analysis involving the stability condition for triples
  we obtain a bound on the values of $\alpha$:  for $n_1 > n_2$ the
  moduli space $\mathcal{M}_{\alpha}$ is empty, unless
  \begin{displaymath}
    0 \leq \alpha \leq \alpha_M =
    \frac{2n_1}{n_1 - n_2}\bigl(\mu(E_1) - \mu(E_2)\bigr).
  \end{displaymath}
  For $n_1 = n_2$ there is no upper bound on $\alpha$, however, the
  moduli spaces $\alpha$ stabilizes for $\alpha$ sufficiently large.
  Thus it makes sense to consider the ``large $\alpha$ moduli space'',
  $\mathcal{M}_{\infty}$, in both cases.  Analysis of
  $\mathcal{M}_{\infty}$ shows that it is non-empty and
  irreducible---for $n_1=n_2$ this is the main result of Markman and
  Xia \cite{markman-xia:2000}.
  
  There is a finite number of so-called critical values of the
  parameter $\alpha$, these are values for which strict
  $\alpha$-semi-stability is possible.  The $\alpha$-stability
  condition remains the same between critical values.  Thus we need to
  study how $\mathcal{M}_{\alpha}$ varies as $\alpha$ crosses a
  critical value and, in particular, show that it remains non-empty
  and irreducible.  The locus where $\mathcal{M}_{\alpha}$ changes
  consists of triples which are strictly $\alpha$-semi-stable for the
  critical value of $\alpha$, and what we need to show is that this
  locus has strictly positive codimension.
  
  The category of triples is an Abelian category and we study strictly
  semi-stable triples via their Jordan-H\"older filtration by stable
  triples.  This involves developing the theory of extensions of
  triples; in particular we show that the set of extensions of a
  triple $T''=(E_{1}'',E_{2}'',\phi'')$ by a triple
  $T'=(E_{1}',E_{2}',\phi')$ is isomorphic to the first
  hypercohomology of the complex
  \begin{align*}
    {E_{1}''}^{*} \otimes E_{1}' \oplus {E_{2}''}^{*} \otimes E_{2}'
    &\to {E_{2}''}^{*} \otimes E_{1}' \\
    (\psi_{1},\psi_{2}) &\mapsto \phi'\psi_{2} - \psi_{1}\phi''.
  \end{align*}
  One can show that the codimension of the locus where
  $\mathcal{M}_{\alpha}$ changes as $\alpha$ crosses a critical value
  is strictly positive if this first hypercohomology group is
  non-vanishing for any extension.  An important fact which we prove
  in the course of these arguments is that for $\alpha \geq 2g-2$ any
  $\alpha$-stable triple is a smooth point of the moduli space.
  
  Finally the non-vanishing of the above first hypercohomology is
  proved using a vanishing criterion which is reminiscent of the one
  described in the proof of Lemma \ref{lemma:minima} above.
\end{proof}

Detailed proofs of these results will appear elsewhere.


\subparagraph{Acknowledgements.}  We thank the mathematics departments
of the University of Illinois at Urbana-Champaign and the Universidad
Aut{\'o}noma de Madrid, the Department of Pure Mathematics of the
University of Porto and the Mathematical Institute of the University
of Oxford for their hospitality during various stages of this
research.  The authors are members of VBAC (Vector Bundles on
Algebraic Curves), which is partially supported by EAGER (EC FP5
Contract no.\ HPRN-CT-2000-00099) and by EDGE (EC FP5 Contract no.\ 
HPRN-CT-2000-00101).  The third author was partially supported by the
Funda{\c c}{\~a}o para a Ci{\^e}ncia e a Tecnologia (Portugal) through
the Centro de Matem{\'a}tica da Universidade do Porto and through
grant no.\ SFRH/BPD/1606/2000.


\end{document}